\documentclass[11pt]{amsart}
\usepackage[dvips]{graphicx}
\usepackage{epsfig}
\usepackage{amscd, amssymb}
\usepackage{times}

\newtheorem{thm}{Theorem}[section]
\newtheorem{lemma}[thm]{Lemma}
\newtheorem{cor}[thm]{Corollary}
\newtheorem{prop}[thm]{Proposition}
\newtheorem{question}[thm]{Question}

\theoremstyle{definition}
\newtheorem{definition}[thm]{Definition}

\newtheorem{remark}[thm]{Remark}

\renewcommand{\(}{\left(}
\renewcommand{\)}{\right)}

\newcommand{\R}{{\mathbf{R}}}		% Real numbers
\newcommand{\C}{{\mathbf{C}}}		% Complex numbers
\newcommand{\Z}{{\mathbf{Z}}}		% Integers
\newcommand{\cn}{\colon}		% colon for $f\colon A \to B$
\newcommand{\la}{\langle}		% Left Angle
\newcommand{\ra}{\rangle}		% Right Angle
\newcommand{\cd}{,\dots,}		% Commas with Dots
\newcommand{\Luko}{L{\" u}k{\H o}}      % Luko Gabor
\newcommand{\Max}{\operatorname{Max}}   % Maximizer set

\DeclareMathOperator{\distort}{distort}

\newcommand{\eq}{~\eqref}
\numberwithin{equation}{section}

\catcode`@=11 \@mparswitchfalse  %This puts the \mnote's on the
\newcounter{mnotecount}[page]

\def\incgr#1#2{\includegraphics[height=#2in]{#1}}

\def\figh#1#2#3%
{ \begin{figure}[ht]
  \begin{center} \incgr{#1}{#2} \end{center}
  \caption{#3} \label{fig:#1}
\end{figure} }

\begin{document}

\title{Circles Minimize most Knot Energies}

\author[Abrams]{Aaron Abrams}
\address{Department of Mathematics, University of Georgia,
Athens, GA 30602}
\email{abrams@math.uga.edu}
\urladdr{www.math.uga.edu/$\sim$abrams/}

\author[Cantarella]{Jason Cantarella${}^1$}
\address{Department of Mathematics, University of Georgia,
Athens, GA 30602}
\email{cantarel@math.uga.edu}
\urladdr{www.math.uga.edu/$\sim$cantarel/}
\thanks{${}^1$ Supported by an NSF Postdoctoral Research Fellowship.}

\author[Fu]{Joseph H. G. Fu${}^2$}
\address{Department of Mathematics, University of Georgia,
Athens, GA 30602}
\email{fu@math.uga.edu}
\urladdr{www.math.uga.edu/$\sim$fu/}
\thanks{${}^2$ Supported by NSF grant DMS-9972094}

\author[Ghomi]{Mohammad Ghomi}
\address{Department of Mathematics, University of South
Carolina, Columbia, SC 29208}
\email{ghomi@math.sc.edu}
\urladdr{www.math.sc.edu/$\sim$ghomi}

\author[Howard]{Ralph Howard${}^3$}
\thanks{${}^3$ Supported in part by DoD Grant No. N00014-97-1-0806}
\address{Department of Mathematics, University of South
Carolina,
Columbia, SC 29208}
\email{howard@math.sc.edu}
\urladdr{www.math.sc.edu/$\sim$howard}

\date{\today}

\begin{abstract}
We define a new class of knot energies (known as {\em renormalization
energies}) and prove that a broad class of these energies are uniquely
minimized by the round circle. Most of O'Hara's knot energies belong
to this class. This proves two conjectures of O'Hara and of Freedman,
He, and Wang. We also find energies not minimized by a round circle.
The proof is based on a theorem of \Luko\ on average chord lengths of
closed curves.
\end{abstract}

\maketitle
%\tableofcontents

%%%%%%%%%%%%%%%%%%%%%%%%%%%%%%%%%%%%% 
%%%%%%%%%%%%%%%%%%%%%%%%%%%%%%%%% 
\section{Introduction} 
%%%%%%%%%%%%%%%%%%%%%%%%%%%%%%%%%%%%% 
%%%%%%%%%%%%%%%%%%%%%%%%%%%%%%%%%
For the past decade, there has been a great deal of interest in
defining new knot invariants by minimizing various functionals on the
space of curves of a given knot type.  For example, imagine a loop of
string bearing a uniformly distributed electric charge, floating in
space. The loop will repel itself, and settle into some least energy
configuration. If the loop is knotted, the potential energy of this
configuration will provide a measure of the complexity of the knot.

In 1991 Jun O'Hara began to formalize this picture 
~\cite{ohara:energy1,ohara:energy2a} by proposing a family of energy
functionals $e_j^p$ (for $j$, $p > 0$) which are based on the
physicists' concept of renormalization, and which are defined by 
$e_j^p[c] := (1/j) (E_j^p[c])^{1/p}$, where
\begin{equation} \label{eq:E}
E_j^p[c] := \iint \left( \frac{1}{|c(s) - c(t)|^j} 
	- \frac{1}{d(s,t)^j} \right)^p \, ds \, dt,
\end{equation}
$c \cn S^1 \to \R^3$ is a unit-speed curve, $|c(s) - c(t)|$ is the
distance between $c(s)$ and $c(t)$ in space, and~$d(s,t)$ is the
shortest distance between $c(s)$ and $c(t)$ along the curve. O'Hara
showed \cite{ohara:energy3} that these integrals converge if the curve
$c$ is smooth and embedded, $j < 2 + 1/p$, and that a minimizing curve
exists in each isotopy class when~$jp > 2$.

It was then natural to try to find examples of these energy-minimizing
curves in various knot types. O'Hara conjectured \cite{ohara:energy2}
in 1992 that the energy-minimizing unknot would be the round circle
for all $e_j^p$ energies with $p \geq 2/j \geq 1$, and wondered
whether this minimum would be unique. Later that year, he provided
some evidence to support this conjecture by proving
\cite{ohara:energy2a} that the limit of $e_j^p$ as $p \to \infty$ and
$j \to 0$ was the logarithm of Gromov's {\em distortion}, which was
known to be minimized by the round circle (see
\cite{kusner&sullivan:distortion} for a simple proof).

Two years later, Freedman, He, and Wang investigated a family of
energies almost identical to the $e_j^p$ energies, proving that the
$e_2^1$ energy was M\"obius-invariant \cite{freedman&he&wang:mobius},
and as a corollary that the overall minimizer for $e_2^1$ was 
the round circle. For the remaining $e_j^1$ energies, they were able
to show only that the minimizing curves must be convex and planar for
$0 < j <3$ (Theorem 8.4). They conjectured that these minimizers were 
actually circles.

We generalize the energies of O'Hara and Freedman-He-Wang as follows:

\begin{definition}
Given a curve $c$ parametrized by arclength, let $|c(s) - c(t)|$ be
the distance between $c(s)$ and $c(t)$ in space, and $d(s,t)$ denote
the shortest distance between $s$ and $t$ along the curve. Given a 
function $F:\R^2 \rightarrow \R$, the energy functional in the form
\begin{equation}
f[c] := \iint F\( |c(s)-c(t)|, d(s,t) \) \, ds \, dt,
\end{equation}
is called the {\em renormalization energy} based on $F$ if it converges
for all embedded $C^{1,1}$ curves.
\end{definition}

The main result of this paper is that a broad class of these energies
are uniquely minimized by the round circle. 

\begin{thm}\label{main}
Suppose $F(x,y)$ is a function from $\R^2$ to $\R$. If $F(\sqrt{x},y)$ is
convex and decreasing in $x$ for $x \in (0,y^2)$ and $y \in (0,\pi)$ 
then the renormalization energy
%\begin{equation*}
%f[c] := \iint F\( |c(s) - c(t)|,d(t,s) \) \, dt\, ds,
%\end{equation*}
based on $F$ is uniquely minimized among closed unit-speed curves of 
length $2\pi$ by the round unit circle.
\end{thm}

It is easy to check that the hypotheses of Theorem~\ref{main} are
slightly weaker than requiring that $F$ be convex and decreasing
in~$x$. The theorem encompasses both O'Hara's and Freedman, He, and
Wang's conjectures:

\begin{cor}\label{E}
Suppose $0 < j < 2 + 1/p$, while $p \geq 1$. Then for every closed
unit-speed curve $c$ in~$\R^n$ with length $2\pi$,
\begin{equation}
E_j^p[c] \geq 2^{3-jp} \pi \int_0^{\frac{\pi}{2}} 
	\left( \(\frac{1}{\sin s}\)^j - \(\frac{1}{s}\)^j \right)^p \, ds.
\end{equation}
with equality if and only if $c$ is the  circle.
\end{cor}

We must include the condition $j < 2 + 1/p$ in our theorem, for otherwise 
the integral defining~$E_j^p$ does not converge. We do not know whether
the condition $p \geq 1$ is sharp, since the energies are well-defined
for $0 < p < 1$, but it is required for our proof. 

We use several ideas from a prophetic paper of \Luko\ G\'abor
\cite{Gabor}, written almost thirty years before the conjectures of O'Hara
and Freedman, He, and Wang were made.  \Luko\footnote{There are references
in the literature to papers authored both by \Luko\ G\'abor and by
G\'abor \Luko.  We are informed that these people are identical and
that \Luko\ is the family name; the confusion likely results from the
Hungarian convention of placing the family name first.} showed that
among closed, unit-speed planar curves of length $2\pi$, circles are
the only maximizers of any functional in the form
\begin{equation}\label{eqn:gabor}
\iint f(|c(s) - c(t)|^2) \, ds \, dt,
\end{equation}
where $f$ is increasing and concave.

Our arguments are modeled in part on Hurwitz's proof of the planar
isoperimetric
inequality~\cite{hurwitz:iso}~\cite[p.\,111]{chern:book}. In
Section~\ref{sec:wint}, we derive a Wirtinger-type inequality
(Theorem~\ref{kostya}), which we use in Section~\ref{sec:dist} to
generalize \Luko's theorem (Theorem~\ref{fchord}).  We then apply this
result to obtain sharp integral inequalities for average chord lengths
and distortions. In the process, we find another proof that the curve
of minimum distortion is a circle. In Section~\ref{sec:energy}, we
give the proof of the main theorem.

All our methods depend on the concavity of $f$ in functionals of the form of 
Equation~\ref{eqn:gabor}. In Section~\ref{sec:expts}, we consider the case where
$f$ is convex, as in the case of the functional
\begin{equation}
\iint |c(s)-c(t)|^p \, ds \, dt
\end{equation}
for $p > 2$. Numerical experiments suggest that the maximizing curve for this 
functional remains a circle for $p < \alpha$, with $3.3 < \alpha < 3.5721$, 
while for $p > 3.5721$, the maximizers form a family of stretched ovals
converging to a doubly-covered line segment as~$p \to \infty$.

%%%%%%%%%%%%%%%%%%%%%%%%%%%%%%%%%%%%%
%%%%%%%%%%%%%%%%%%%%%%%%%%%%%%%%%
\section{A  Wirtinger type inequality}
\label{sec:wint}
%%%%%%%%%%%%%%%%%%%%%%%%%%%%%%%%%%%%%
%%%%%%%%%%%%%%%%%%%%%%%%%%%%%%%%%

\begin{definition}\label{def:lambda}
Let $\lambda\cn \R\to \R$ be given by
\begin{equation}\label{eq:lambda}
\lambda(s) := 2 \sin \frac{s}{2}.
\end{equation}
For $0\leq s\leq 2\pi$, $\lambda(s)$ is the length of the chord
connecting the end points of an arc of length $s$ in the unit circle.
\end{definition}

Our main aim in this section is to prove the following inequality,
modeled after a well known lemma of
Wirtinger~\cite[p.~111]{chern:book}.  For simplicity, we restrict our
attention to closed curves of length~$2\pi$ in $\R^n$.

% We give two
%proofs, one based on the elementary theory of Fourier series, and one
%based on the maximum principle for ordinary differential equations.
%
%Our main aim in this section is to prove the following inequality, 
%modeled after a well known lemma of Wirtinger~\cite[p.\,111]{chern:book}.
%For simplicity, we will restrict our attention to closed curves 
%of length $2\pi$ in $\R^n$.

\begin{thm}\label{kostya}
Let $c\cn S^1:= \R/2\pi \Z \to \R^n$ be an absolutely continuous function. 
If $c'(t)$ is square integrable, then for any $s\in\R$
\begin{equation}\label{wint-ineq}
\int |c(t+s)-c(t)|^2\,dt
\le \lambda^2(s) \int |c'(t)|^2\,dt,
\end{equation}
with equality  if and only if $s$ is an integral multiple of $2\pi$ or
\begin{equation}\label{ellipse}
c(t)=a_0+(\cos t)\,a+(\sin t)\,b
\end{equation}
for some vectors $a_0,a,b\in \R^n$.
\end{thm}

We give two proofs of this result, one based on the elementary theory
of Fourier series, and one based on the maximum principle for ordinary
differential equations. 

\begin{proof}[Fourier series proof]
We assume that $c\cn S^1 \to \R^n\subset \C^n$, as the complex form of
the Fourier series is more convenient.  $\C^n$ is equipped with its
standard positive definite Hermitian inner product $\la
v,w\ra=\sum_{k=1}^nz_k\overline{w}_k$ where $v=(v_1\cd v_n)$ and
$w=(w_1\cd w_n)$. This agrees with the usual inner product on
$\R^n\subset\C^n$.  The norm of $v\in \C^n$ is given by
$|v|:=\sqrt{\la v,v\ra}$, and $i:=\sqrt{-1}$.

The facts about Fourier series required for the proof are as
follows. If $\phi\cn S^1\to\C^n$ is locally square integrable then it
has a Fourier expansion
$$
\phi(t)=\sum_{k=-\infty}^\infty \phi_k e^{kti},
$$
(the convergence is in $L^2$ and the series may not converge pointwise). 
The $L^2$ norm of $\phi$ is given by
\begin{equation}\label{eq:fourier1}
\int |\phi(t)|^2\,dt
	= 2\pi \sum_{k=-\infty}^\infty |\phi_k|^2.
\end{equation}
If $\phi$ is absolutely continuous and $\phi'$ is locally square
integrable then $\phi'$ has the Fourier expansion $\phi'(t)=
i \sum_{k=-\infty}^\infty k\phi_k e^{kti}$ and therefore
\begin{equation}\label{eq:fourier2}
\int |\phi'(t)|^2\,dt= 2\pi \sum_{k=-\infty}^\infty k^2 |\phi_k|^2
= 2\pi \sum_{k=1}^\infty k^2(|\phi_{-k}|^2+|\phi_k|^2),
\end{equation}
as the contribution to the middle sum from the term $k=0$ is zero.

Let $\sum_{k=-\infty}^\infty a_ke^{kti}$ be the Fourier
expansion of $c(t)$, where $a_k\in \C^n$.  Then
\begin{align*}
c(t+s/2)-c(t-s/2)
&=\sum_{k=-\infty}^\infty\(e^{ksi/2}-e^{-ksi/2}\)a_ke^{kti}\\
&=2i\sum_{k=-\infty}^\infty \( \sin \frac{ks}{2} \) a_ke^{kti}.
\end{align*}
Therefore, using \eqref{eq:fourier1}, we have
\begin{align}\label{eq:int1}
\int |c(t+s)-c(t)|^2\,dt
&= \int \left|c\(t + \frac{s}{2}\) - c\(t - \frac{s}{2}\)\right|^2 \, dt \nonumber \\
& =2\pi |2i|^2 \sum_{k=-\infty}^\infty \( \sin^2 \frac{ks}{2} \) |a_k|^2 \nonumber\\
&=8\pi \sum_{k=1}^\infty \( \sin^2 \frac{ks}{2} \) \( |a_{-k}|^2+|a_k|^2 \).
\end{align}
Also, by \eqref{eq:fourier2} and \eqref{eq:lambda},
\begin{align}\label{eq:int2}
\lambda^2(s) \int |c'(t)|^2\,dt
	&=\(4 \sin^2 \frac{s}{2} \) 
	  \( 2\pi \sum_{k=1}^\infty k^2 (|a_k|^2+|a_{-k}|^2) \) \nonumber\\
        &=8\pi \sum_{k=1}^\infty \( k^2 \sin^2 \frac{s}{2} \) (|a_k|^2+|a_{-k}|^2).
\end{align}
Subtracting \eqref{eq:int1} from \eqref{eq:int2}, we set
\begin{align*}
\rho_c(s)&:=\lambda^2(s)\int |c'(t)|^2\,dt
	-\int |c(t+s)-c(t)|^2\,dt\\
         &=8\pi \sum_{k=2}^\infty \( k^2 \sin^2 \frac{s}{2} - \sin^2 \frac{ks}{2} \)
	   (|a_{-k}|^2+|a_k|^2).
\end{align*}

Lemma \ref{trig-ineq} (below) implies that
$\rho_c(s) \ge 0$ with equality if and only if $s$ is a multiple of $2\pi$, or
$a_k=a_{-k}=0$ for all $k\ge 2$. The latter occurs if and only if
\begin{equation}\label{eq:c}
c(t)=a_{-1}e^{-it}+a_0+a_1 e^{it}
    =a_0+(\cos t)\,a+(\sin t)\,b
\end{equation}
where $a:=a_{1}+a_{-1}$ and $b:=i(a_{1}-a_{-1})$.
\end{proof}

\begin{lemma}\label{trig-ineq}
Let $k\ge 2$ be an integer. Then
\begin{equation}\label{eq:theta}
\sin^2(k\theta)\leq k^2\sin^2(\theta),
\end{equation}
with equality if and only if $\theta=m\pi$ for some integer $m$
\end{lemma}

\begin{proof} 
If $\theta=m\pi$, for some integer $m$, then equality holds
in~\eqref{eq:theta}. If $\theta$ is not an integer multiple of~$\pi$, we set
$q_k(\theta):=|\sin(k\theta)/\sin(\theta)|$. Then $|\cos(\theta)|<1$,
and the addition formula for sine yields
\begin{equation}\label{eq:q}
q_{k+1}(\theta) 
 =|\cos(\theta)\,q_k(\theta)+\cos(k\theta)|
 <  q_k(\theta)+1,
\end{equation}
Since $q_1(\theta)\equiv1$, we then have $q_k(\theta) < k$ by induction, which completes the proof.
\end{proof}

\begin{proof}[Maximum principle proof]
This method is an adaptation of \Luko's original
approach~\cite{Gabor}.  In that paper, he solves a discrete version of
the problem, showing that the average squared distance between the
vertices of an $n$-gon of constant side length is maximized by the
regular $n$-gon.  He then obtains the main result by approximation. We go
directly to the continuum case, which turns out to be simpler.

To simplify notation, let $L=\int |c'(t)|^2 \,dt$.  Let
\begin{gather*}
f(s):=\int |c(t+s)-c(t)|^2 \,dt,\\
\Lambda(s):=\lambda^2(s)\int |c'(t)|^2 =L\lambda^2(s).
\end{gather*}
We  claim that $f$ is $C^2$ with
\begin{gather*}
f'(s)=2 \int \la c(t)-c(t-s), c'(t)\ra \,dt,\\
f''(s)=2\int \la c'(t-s), c'(t)\ra \,dt,
\end{gather*}
and initial conditions
\begin{equation}\label{0.3}
f(0)=0,\qquad f'(0)=0,\qquad f''(0)=2\int|c'(t)|^2\,dt= 2L.
\end{equation}
These formulas are clear when  $c$ is $C^2$ and hold in
the general case by approximating by $C^2$ functions.  The explicit
formula for $f''$ makes it clear that $f$ is $C^2$.

Next we  derive a differential inequality for $f$, using an
elementary geometric fact (which appears in a slightly different form in 
\Luko's paper as Lemma~7):

\begin{lemma} \label{lem:luko}
For any tetrahedron $A$, $B$, $C$, $D$ in $\R^n$,
\begin{equation} \label{eqn:luko}
|AC|^2 + |BD|^2 \leq |BC|^2 + |AD|^2 + 2 |AB|\, |CD|,
\end{equation}
with equality if and only if $AB$ and $DC$ are parallel as
vectors.
\end{lemma}

\begin{proof}
Denote the vectors $AB$, $BC$, $CD$, $DA$ by $v_1$, $v_2$, $v_3$, $v_4$. Then
$\sum v_i = 0$, and 
\begin{align*}
|AC|^2 + |BD|^2 &= \frac{1}{2}\left( |v_1 + v_2|^2 + |v_2 + v_3|^2 + |v_3 + v_4|^2 + |v_4 + v_1|^2 \) \\
		&= \sum_{i=1}^{4} |v_i|^2 
			+ \la v_1,v_2 \ra + \la v_2,v_3 \ra + \la v_3,v_4 \ra + \la v_4,v_1 \ra \\
		&= \sum_{i=1}^{4} |v_i|^2 + \la v_1 + v_3, v_2 + v_4 \ra \\
		&= \sum_{i=1}^{4} |v_i|^2 - |v_1 + v_3|^2 \\
		&\leq \sum_{i=1}^{4} |v_i|^2 - (|v_1| - |v_3|)^2 \\
		&= |v_2|^2 + |v_4|^2 + 2 |v_1| |v_3| = |BC|^2 + |AD|^2 + 2 |AB|\, |CD|. 
\end{align*}
Equality holds if and only if $v_3=-\rho v_1$ for some $\rho>0$, which
is equivalent to $AB$ and $DC$ being parallel as vectors.
\end{proof}
For any $t$, $s$ and $h$, we can apply Lemma~\ref{lem:luko} to the
tetrahedron $c(t)$, $c(t+s+h)$, $c(t+s)$, $c(t+h)$ to derive the
equation
\begin{align*}
|c(t+s)-c(t)|^2 & + |c(t+s+h)-c(t+h)|^2  \\
&\le  |c(t+s+h)-c(t+s)|^2 + |c(t+h)-c(t)|^2  \\
&\quad + 2|c(t+s +h) -c(t)| \, |c(t+s) -c(t+h)|.  
\end{align*}
Holding $s, h$ fixed and integrating with respect to $t$,
\begin{align*}
2 f(s) &\le 2 f(h) + 2 \int |c(t+s+h)-c(t)| \, |c(t+s) -c(t+h)| \, dt \\
       &\le 2 f(h) + 2 \sqrt{ f(s+h) f(s-h)} \label{0.4}
\end{align*}
by the Cauchy-Schwartz inequality. Therefore $f(s) \le f(h)+ \sqrt{f(s+h )
f(s-h)}$. For any fixed~$s$, this can be rewritten
\begin{equation*}
g(h) := \frac{1}{2} \big(\log f(s+h) + \log f(s-h)\big) - \log\big(f(s) - f(h)\big)
\geq 0.
\end{equation*}
When $s$ is not a multiple of $2\pi$, $f(s) > 0$ and $g$ is well-defined for
small $h$. Further, $g$ has a local minimum at $h=0$, and so the second derivative
of $g$ is non-negative at zero. Using~\eqref{0.3}, this tells us that
\begin{equation}\label{0.5}
\frac{d^2}{ds^2}\log f(s) \ge \frac{-2L}{f(s)}.
\end{equation}
Meanwhile, $\Lambda(s)$ satisfies the differential equation
\begin{equation}\label{0.6}
\frac{d^2}{ds^2}\log \Lambda(s) = \frac{-2L}{\Lambda(s)}.
\end{equation}

We are trying to show that $f(s) \le \Lambda(s)$ and that if equality holds 
for any $s\in (0,2\pi)$, then $f(s) \equiv \Lambda(s)$. Let 
\begin{equation*}
u(s)=\log\frac{f(s)}{\Lambda(s)}=\log f(s)-\log \Lambda(s).
\end{equation*}
In these terms, we want to show that $u(s) \le 0$ and that if $u(s)=0$ for some 
$s \in (0,2\pi)$ then $u \equiv 0$.  Using~\eqref{0.5} and~\eqref{0.6},
\begin{equation*}
u''(s) \ge \frac{-2L}{f(s)}+\frac{2L}{\Lambda(s)}
	=\frac{2L}{f(s)}\(\frac{f(s)}{\Lambda(s)}-1\)
	=\frac{2L}{f(s)}\(e^{u(s)}-1\)
	\ge \frac{2L}{f(s)}u(s).
\end{equation*}
By two applications of L'Hospital's rule, we compute $\lim_{s\to0} u(s) = 0$.
Thus $\lim_{s\to2\pi} u(s) = 0$, as well. So if $u$ is ever positive,
it will have a positive local maximum at some point $s_0 \in (0,2\pi)$. At that
point,
\begin{equation*}
0 \geq u''(s_0) \geq \frac{2L}{f(s_0)} u(s_0) > 0,
\end{equation*}
which is a contradiction. So $u$ is non-positive on $(0,2\pi)$.
Further, if $u$ is zero at any point in $(0,2\pi)$, the strong maximum
principle~\cite[Thm~17 p.\,183]{Spivak:comp5} implies that $u$ vanishes 
on the entire interval. Thus $f(s) \leq \Lambda(s)$ with equality
at any point of $(0,2\pi)$ if and only if $f(s) \equiv \Lambda(s)$.

Last, we show that if $f(s) = \int |c(t+s)-c(t)|^2 \,dt \equiv
\lambda^2(s) \int |c'(t)|^2 \,dt = \Lambda(s),$ then $c$ is an 
ellipse. By our work
above, if $f = \Lambda$, then for each fixed $s$, $c$ maximizes 
$\int |c(t+s)-c(t)|^2 \,dt$ 
subject to the constraint that $\int |c'(t)|^2 \,dt$ is held
constant.  The Lagrange multiplier equation for this variational
problem is
\begin{equation*}
c''(t)=M\big(c(t+s)-2c(t)+c(t-s)\big)
\end{equation*}
where $M$ is a constant depending on $s$.  When $s=\pi$ we can use
the fact that $c$ has period~$2\pi$ and this becomes
\begin{equation*}
c''(t)=2M\big(c(t+\pi)-c(t)\big).
\end{equation*}
Differentiating twice with respect to $t$, and using both the periodicity and the equation,
\begin{align*}
c''''(t)&= 2M   \big(c''(t+\pi)-c''(t)\big)\\
	&= 4M^2 \big(c(t)-c(t-\pi)-c(t+\pi)+c(t)\big)\\
	&=-8M^2 \big(c(t+\pi)-c(t)\big)\\
	&=-4Mc''(t).
\end{align*}
So $c''$ satisfies the equation $g'' = -4M g$ and has period $2\pi$.
This implies that $4M = k^2$ for some $k\in\Z$, and 
$c''(t) = (\cos kt) V + (\sin kt) W$ with 
$V$ and $W$ in $\R^n$.  But $k=\pm 1$, for otherwise
$f(2\pi/k)=0\ne \Lambda(2\pi/k)$, a contradiction.  Taking two 
antiderivatives,
\begin{equation}
c(t) = a_0 + tb_0 + \(\cos t\) \, a + \(\sin t\) \, b,
\end{equation}
with $a_0, b_0, a, b$ in $\R^n$. Periodicity implies that $b_0=0$, completing
the proof.
\end{proof}

\begin{remark}\label{affine}
By equation \eqref{eq:c}, extremals for the inequality of
Theorem~\ref{kostya} are either ellipses or double coverings of line
segments, depending on whether $a$ and $b$ are linearly independent.
Thus the set of extremal curves is invariant under affine maps of
$\R^n$. When the extremal is an ellipse, the parameterization is a
constant multiple of the \emph{special affine arclength}
(c.f. \cite[p.\,7]{Burago-Zal:geometric-ineq},
\cite[p.\,56]{Spivak:comp2}).  It would be interesting to find an
affine invariant interpretation of inequality\eq{wint-ineq} or of the
deficit $\rho_c(s)$ used in the first proof---especially when $c$ is a
convex planar curve.
\end{remark}

%%%%%%%%%%%%%%%%%%%%%%%%%%%%%%%%%%%%%
%%%%%%%%%%%%%%%%%%%%%%%%%%%%%%%%%%%%%
\section{Inequalities for Concave Functionals}
\label{sec:dist}
%%%%%%%%%%%%%%%%%%%%%%%%%%%%%%%%%%%%%
%%%%%%%%%%%%%%%%%%%%%%%%%%%%%%%%%%%%%

We now apply Theorem~\ref{kostya} to obtain an inequality for chord
lengths.  Recall Definition~\ref{def:lambda}, that $\lambda(s)$ is the
length of a chord of arclength $s$ on the unit circle.

\begin{thm}\label{fchord}
Let $c$ be a closed, unit-speed curve of length $2\pi$ in $\R^n$.  For
$0<s<2\pi$, if $f\cn\R \to \R$ is increasing and concave on
$(0,d(0,s)^2]$, where $d(s,t)$ is the shortest distance along the
curve between $c(s)$ and $c(t)$, then
\begin{equation}\label{ineq:chord}
\frac{1}{2\pi}\int f\(|c(t+s)-c(t)|^2\)\,dt
  \le f\(\lambda^2(s)\)
\end{equation}
and equality holds if and only if $c$ is the unit circle.

%Note: The last half of this theorem is now encompassed by our main
%theorem. So we simply defer the proof until the next section.

%Further, if $f\cn\R \to \R$ is increasing and concave on $(0,\pi^2]$, 
%\begin{equation}\label{ave:chord}
%\frac{1}{4\pi^2}\iint f\(|c(t)-c(s)|^2\)\,dt\,ds
%  \le \frac{1}{2\pi}\int f\(\lambda^2(s)\)\,ds.
%\end{equation}
%If $ \frac{1}{2\pi}\int f\(\lambda^2(s)\)\,ds$ is finite, then
%equality holds if and only if $c$ is the unit circle.
\end{thm}

\begin{proof}
The shortest distance between $c(t)$ and $c(t+s)$ along the curve is $d(0,s)$. Thus,
the squared chord length $|c(t+s) - c(t)|^2$ is in $(0,d(0,s)^2]$, except when $s=0$.
Being undefined at this point does not affect the existence of the integrals.
Using Jensen's inequality for concave functions~\cite[p.~115]{Royden:analysis},
Theorem~\ref{kostya}, that $f$ is increasing, and that $|c'(t)|=1$ for
almost all~$t$, we have
\begin{align*}
\frac1{2\pi}\int f\(|c(t+s)-c(t)|^2\)\,dt
	&\le f\(\frac{1}{2\pi} \int |c(t+s)-c(t)|^2\,dt \)\\
	&\le f\(\frac{\lambda^2(s)}{2\pi}\int |c'(t)|^2 \,dt\)\\
	&=f\(\lambda^2(s)\).
\end{align*}
If equality holds in \eqref{ineq:chord}, then the above string of
inequalities implies that equality holds between the two middle terms,
i.e., equality holds in \eqref{wint-ineq}.  Thus, since $0 < s < 2\pi$, we
may apply Theorem~\ref{kostya} to conclude that $c(t)$ must be as in
\eqref{ellipse}.  Since $c$ has unit speed, it follows that
\begin{equation*}
c'(t)= -\(\sin t \)\, a + \(\cos t\) \,b
\end{equation*}
is a unit vector for all $t$, which forces the vectors $a$ and $b$ to
be orthonormal, and so implies that $c$ is the unit circle.
Conversely, if $c$ is the unit circle, then $|c(t+s)-c(t)|=\lambda(s)$
for all $t$ and therefore equality holds in\eq{ineq:chord}.

%Note: Again, the last part of theorem will be taken care of in 
%the main theorem in the next section.

%It is clear that $(0,d(0,s)^2] \subset (0,\pi^2]$ for
%all $s \in (0,2\pi)$. The inequality \eqref{ave:chord} now follows
%easily by a change of variables and applying \eqref{ineq:chord}:
%\begin{align*}
%\frac{1}{4\pi^2} \iint f\( |c(t)-c(s)|^2 \)\,dt\,ds
%	&=\frac1{4\pi^2}\iint f\( |c(t+s)-c(t)|^2 \)\,dt\,ds\\
%	&\le \frac{1}{2\pi} \int f\(\lambda^2(s)\)\,ds.
%\end{align*}
%This shows that if the last of these integrals is finite, then equality
%holds in~\eqref{ave:chord} if and only if equality holds in
%\eqref{ineq:chord}, which can only happen when $c$ is the unit circle.
\end{proof}

Letting $f(x)=\sqrt{x}$ in Theorem~\ref{fchord}, we obtain the
following inequality:

\begin{cor}\label{chord}
Let $c$ be a closed, unit-speed curve of length $2\pi$ in $\R^n$. Then for
any $s\in (0,2\pi)$,
\begin{equation}
\frac1{2\pi}\int |c(t+s)-c(t)| \,dt \le \lambda(s),\label{pre-dist}
\end{equation}
with equality if and only if $c$ is the unit circle.\qed
\end{cor}

Next we apply Theorem~\ref{fchord} to obtain sharp inequalities for
Gromov's \emph{distortion}
\cite{gromov:filling,kusner&sullivan:distortion}.  By definition, the
distortion of a curve is the maximum value of the ratio of the
distance in space to the distance along the curve for all pairs of
points on the curve. As we mentioned above, distortion is a limit of
O'Hara energies: $\exp(e_0^\infty(c))=\distort(c)$
\cite[p.\,150]{ohara:energy3}.

The inequality \eqref{eq:distortion2} is due to Gromov
\cite[pp.\,11--12]{Gromov:metric}, \cite{kusner&sullivan:distortion}.
As always, while we state our results for curves of length $2\pi$, the
corresponding result holds for curves of arbitrary length.

\begin{cor}\label{s-dist}
For every closed, unit-speed curve $c$ of length $2\pi$ in $\R^n$
\begin{align}
\distort_s(c):=\sup_{t\in \R} \frac{s}{|c(t+s)-c(t)|}& \ge\frac{s}{\lambda(s)},\label{eq:distortion}\\
\distort(c)  :=\sup_{s\in (0,\pi]} \sup_{t\in\R}
\frac{s}{|c(t+s)-c(t)|}&\ge\frac{\pi}{2}\label{eq:distortion2},
\end{align}
with equalities if and only if $c$ is the unit circle.  
\end{cor}
\begin{proof}
In both cases equality is clear for the unit circle.
By the mean value property of integrals and inequality\eq{pre-dist},
\begin{align*}
\frac{1}{\distort_s(c)}&=\inf_{t\in \R}\frac{|c(t+s)-c(t)|}{s}\le
\frac{1}{2\pi s}\int |c(t+s)-c(t)|\,dt\le \frac{\lambda(s)}{s},
\end{align*}
establishing~\eqref{eq:distortion}.  Further, 
equality in~\eqref{eq:distortion} implies equality in~\eqref{pre-dist}, 
which, by Theorem~\ref{fchord}, 
happens if and only if $c$ is the unit circle. 

The proof of~\eqref{eq:distortion2} follows easily from
\eqref{eq:distortion}:  
\begin{equation*}
\distort(c)=\sup_{s\in (0,\pi]} \distort_s(c) \ge
\distort_{\pi}(c)\ge\frac{\pi}{\lambda(\pi)}=\frac{\pi}{2},
\end{equation*}
and again equality implies in particular that
$\distort_{\pi}(c)=\pi/\lambda(\pi)$, which, by \eqref{eq:distortion},  
happens if and only if $c$ is the unit circle.
\end{proof}

For general maps $f\cn M\to\R^n$ of a compact Riemannian manifold to
Euclidean space Gromov~\cite[p.\,115]{gromov:filling} has given, by
methods related to ours, lower bounds---which are not sharp---for the
distortion of $f$ in terms of the first eigenvalue of $M$ and the
average square distance, $\operatorname{Vol}(M)^{-2}\iint_{M\times M}
d(x,y)^2\,dx\,dy$, between points of $M$ (where $d$ is the Riemannian
distance).

%%%%%%%%%%%%%%%%%%%%%%%%%%%%%%%%%%%%%
%%%%%%%%%%%%%%%%%%%%%%%%%%%%%%%%%%%%%
\section{Proof of the Inequality for  Energies}
\label{sec:energy}
%%%%%%%%%%%%%%%%%%%%%%%%%%%%%%%%%%%%%
%%%%%%%%%%%%%%%%%%%%%%%%%%%%%%%%%%%%%

We are now ready to prove the main theorem. We start by restating it.
\begin{thm} \label{thm:main}
Suppose $F(x,y)$ is a function from $\R^2$ to $\R$. If $F(\sqrt{x},y)$ is
convex and decreasing in $x$ for $x \in (0,y^2]$ for all $y \in (0,\pi)$ 
then the renormalization energy based on $F$
\begin{equation*}
 f[c] := \iint F\( |c(s) - c(t)|,d(t,s) \) \, dt\, ds,
\end{equation*}
is uniquely minimized among closed unit-speed curves of length $2\pi$ 
by the round unit circle.
\end{thm}

\begin{proof}
Making the substitution $s \mapsto s - t$, $t \mapsto t$, changing the order of integration,
and using the fact that $d(s,t) = d(s+a,t+a)$ for any $a$, we have
\begin{equation*}
  \iint F\(|c(s) - c(t)|,d(s,t) \) \, ds \, dt = 
  \iint F\(|c(t+s) - c(t)|,d(0,s) \) \, dt \,ds.
\end{equation*}
For each $s \in (0,2\pi)$, if we let $f(x) = -F(\sqrt{x},d(0,s))$, then
\begin{equation*}
  \int F\( |c(t+s) - c(t)|, d(0,s) \) \, dt = -\int f \( |c(t+s) - c(t)|^2 \) \, dt
\end{equation*}
and $f$ is increasing and concave on $(0,d(0,s)^2]$. By Theorem~\ref{fchord},
\begin{equation}\label{eq:last}
  -\int f\(|c(t+s) - c(s)|^2\) \, dt \geq - 2\pi f\( \lambda^2(s) \),
\end{equation}
with equality if and only if $c$ is the unit circle. Integrating this from $s=0$ to
$s=2\pi$ tells us that $f[c]$ is greater than or equal to the corresponding value
for the unit circle, with equality if and only if \eqref{eq:last} holds for almost all
$s \in [0,2\pi]$. But if equality holds for any $s \in (0,2\pi)$, then $c$ is the unit
circle.
\end{proof} 

We now prove the corollary.

\begin{cor}\label{cor:Ohara}
Suppose $0 < j < 2 + 1/p$, while $p \geq 1$. Then for every closed
unit-speed curve $c$ in~$\R^n$ with length $2\pi$,
\begin{equation}\label{eq:Ohara}
E_j^p[c] \geq 2^{3-jp} \pi \int_0^{\frac{\pi}{2}} 
	\left( \(\frac{1}{\sin s}\)^j - \(\frac{1}{s}\)^j \right)^p \, ds.
\end{equation}
with equality if and only if $c$ is the  circle.
\end{cor}

\begin{proof}
If we let 
\begin{equation*}
 F(x,y):= \(\frac{1}{x^{j}}-\frac{1}{y^j}\)^p,
\end{equation*}
then using \eqref{eq:E}, we see that $E_j^p[c]$ is the renormalization
energy based on $F$. We must show that $F(\sqrt{x},y)$ is convex and
decreasing in $x$ for $x \in (0,y^2]$ for all $y \in (0,\pi)$.  It
suffices to check the signs of the first and second partial
derivatives of $F(\sqrt{x},y)$ with respect to $x$ on $(0,y^2)$.

When $p\ge1$, $y \neq 0$, and $x\in(0,y^2)$,
\begin{equation*}
\frac{\partial F(\sqrt{x},y)}{\partial x}=-\frac{jp}{2x^{(j+2)/2}}
      \(\frac{1}{x^{j/2}}-\frac{1}{y^j}\)^{p-1} < 0,
\end{equation*}
and
\begin{equation*}
\frac{\partial^2 F(\sqrt{x},y)}{\partial x^2}=
 \frac{j(j+2)p}{4x^{(j+4)/2}}
	\(\frac{1}{x^{j/2}}-\frac{1}{y^j}\)^{p-1} +
 \frac{j^2p(p-1)}{4x^{(j+2)}}
	\(\frac{1}{x^{j/2}}-\frac{1}{y^j}\)^{p-2} > 0.
\end{equation*}
Since $x^{j/2}$ can be arbitrarily close to $y^j$ if the curve is
nearly straight, examining this equation shows that the condition
$p \ge 1$ is required to enforce the convexity of $F(\sqrt{x},y)$.

So for every $y \neq 0$, $F(\sqrt{x},y)$ is decreasing and convex on
$(0,y^2]$.  Further, a direct calculation shows that $\int
F\(\lambda^2(s),s\)\,ds<\infty$ when $j<2+{1}/{p}$.

Thus $F$ satisfies the hypotheses of Theorem~\ref{thm:main}. Computing
the energy of the round circle by changing the variable $s \mapsto 2s$
and noting that the resulting integrand is symmetric about $s =
\pi/2$, we have
\begin{align*}
E_j^p[c]
&\ge 2\pi \int F\(\lambda^2(s),d(0,s) \)\,ds\\
&=   2^{2-jp}\pi \int_0^\pi \(\(\frac{1}{\sin s}\)^j-
               \(\frac{1}{\min\{s, \pi-s\}}\)^j\)^p\,ds\\
&=  2^{3-jp}\pi \int_0^{\pi/2}\(\(\frac{1}{\sin s}\)^j-\(\frac{1}{s}\)^j\)^p\,ds\end{align*}
with equality if and only if $c$ is the unit circle. 
\end{proof}

\section{Convex functionals and numerical experiments \label{sec:expts}}

All of our work so far has depended on the hypotheses of 
Theorem~\ref{fchord}: our energy integrands must be increasing, {\em
concave} functions of squared chord length. It is this condition which
restricts Corollary~\ref{cor:Ohara} to $e_j^p$ energies with $p \geq 1$.
To investigate the situation where $p < 1$, we focus our attention on 
a model problem. If $0<p<2$, then $f(x)=x^{p/2}$ is increasing and concave; so
Theorem~\ref{thm:main} implies that among closed, unit speed curves
of length $2\pi$ in $\R^n$,
\begin{equation*}
A_p[c]:=\(\frac{1}{4\pi^2}\iint |c(t)-c(s)|^p\,dt\,ds\)^{\frac1p} \le
\(\frac1{2\pi}\int \(\lambda(s)\)^p\,ds\)^{\frac1p},
\end{equation*}
where equality holds if and only if $c$ is the unit circle. When $p = 1$,
this inequality corresponds to the theorem of \Luko~\cite{Gabor} mentioned
in the introduction. It is natural to ask:
\begin{question}
Which closed, unit speed curves of length $2\pi$ maximize $A_p$ for $p > 2$?
\end{question}

We begin by sketching a proof that such a maximizing curve exists for $p>0$.

\begin{prop}\label{prop:exists}
Let $A_p[c]$ be defined as above. For $p>0$, there exists a closed,
unit-speed curve of length $2\pi$ maximizing $A_p[c]$. Further, every
maximizer of $A_p[c]$ is convex and planar.
\end{prop}

\begin{proof}
Sallee's stretching theorem~\cite{sallee} (see
also~\cite{Ghomi-Howard;unfold}) says that for any closed
unit-speed space curve~$c$ of length~$2\pi$, there exists a
corresponding closed, convex, unit-speed plane curve $c^*$ of
length~$2\pi$ such that for every $s$, $t$ in $[0,2\pi]$,
\begin{equation}
	|c(t) - c(s)| \leq |c^*(t) - c^*(t)|,
\end{equation}
with equality for all~$s$ and~$t$ iff~$c$ is convex and planar.  Since
the integrand defining~$A_p[c]$ is an increasing function of chord
length for $p>0$, this implies that every maximizer of~$A_p[c]$ must
be convex and planar. 

Let $\mathcal{U}$ denote the space of closed, convex, planar,
unit-speed curves of length~$2\pi$ which pass through the origin, with
the~$C^0$ norm. It now suffices to show that a maximizer of~$A_p[c]$
exists in $\mathcal{U}$.

Blaschke's selection principle~\cite[p. 50]{schneider:book} implies
that this space of parametrized curves is compact in the~$C^0$
norm. It easy to see that $A_p[c]$ is $C^0$-continuous for $c$ in
$\mathcal{U}$ (in fact, it is jointly continuous in $p$ and $c$ on the
product $(0,\infty)\times \mathcal{U}$),
completing the proof.
\end{proof}

We conjecture that these maximizers are unique (up to rigid motions),
and depend continuously on $p$. It is easy to see the following:
\begin{lemma}
\label{lem:compact} 
As above, let $\mathcal{U}$ denote the space of closed, convex,
planar, unit-speed curves of length~$2\pi$ with the~$C^0$ norm.  Then
\begin{equation*}
\Max := \{(p,c_p) ~|~ \text{$c_p$ is a maximizer of $A_p$} \} \subset (0,\infty)
\times \mathcal{U} 
\end{equation*}
is locally compact and projects onto $(0,\infty)$.
\end{lemma}

\begin{proof}
We know from the proof of Proposition~\ref{prop:exists} that $A$ is a 
$C^0$-continuous functional on the space $(0,\infty) \times
\mathcal{U}$.
%\mnote{\ralph A bit of rewriting following Joe's remarks.}
%More
%generally it is easy to see that the map $(p,c)\mapsto A_p[c]$ is
%continuous on the product $[0,\infty)\times
%\mathcal{U}$.
%\mnote{\ralph The proof seems to require continuity on
%the product and not just continuity in each of the variables
%separately.}
%
%It is easy
%to see that $A_p$ is continuous with respect to $p$, as well.
If we choose any~$(p_0,c_{p_0})$, and choose a compact interval $I \subset
\R$ containing~$p_0$, then $\Max_I = \{(p,c_p) \in \Max~|~p \in I\}$ contains
a neighborhood of~$(p_0,c_{p_0})$. We now show~$\Max_I$  is compact.

Take any sequence~$(p_i,c_{p_i}) \in \Max_I$. Since~$I$ is compact, we
may assume that the~$p_i$ converge to some $p$. Since~$\mathcal{U}$ is
$C^0$-compact (see the proof of Proposition~\ref{prop:exists}), we may
also assume that the~$c_{p_i}$ converge to some~$c$. It remains to show
that~$c$ is a maximizer for~$A_p$.

If not, there exists some $c_p$ with~$A_p[c_p] > A_p[c]$. But then
\begin{equation*}
  \lim_{i \to \infty} A_{p_i}[c_p] = A_p[c_p] > A_p[c] = \lim_{i\to\infty} A_{p_i}[c_{p_i}],
\end{equation*}
since~$A_p$ is continuous in~$p$. On the other hand, since the~$c_{p_i}$ are 
maximizers for the~$A_{p_i}$, we have~$A_{p_i}[c_{p_i}] \geq A_{p_i}[c_p]$ for
each~$i$, and so
\begin{equation*}
  \lim_{i \to \infty} A_{p_i}[c_p] \leq \lim_{i\to\infty} A_{p_i}[c_{p_i}].
\end{equation*}
\end{proof}

Together with uniqueness, this would prove that the set $\Max$ is a
single continuous family of curves depending on $p>0$. As it stands,
Lemma~\ref{lem:compact} tells us surprisingly little about the
structure of $\Max$. For instance, there are locally compact subsets 
of $\R^2$ which project onto the positive $x$-axis but which are 
totally disconnected; one example is
%example, the subset of $\R^2$ defined by
%\mnote{\ralph Following Joe's remark the set has been given as set
%of ordered pairs.}
$$
\left\{\bigg(\sum_{i\ge N} \frac{a_i}{3^i}\, ,   \sum_{\{i\, |\,a_i = 1\}}
\frac{a_i}{3^i} \bigg)\, \bigg|\  a_i \in \{0,1,2\}, N \in \Z \right\} .
$$
%the multi-valued function
%\begin{equation*}
%x = \sum_i \frac{a_i}{3^i},\,\,\, (a_i \in \{0,1,2\}) \quad f(x) = \sum_{\{i |
%a_i = 1\}} \frac{a_i}{3^i},
%\end{equation*}

\begin{figure}[ht]
  \hbox to 6in{\hfill
  \vbox {\hbox{\includegraphics[height=1.5in]{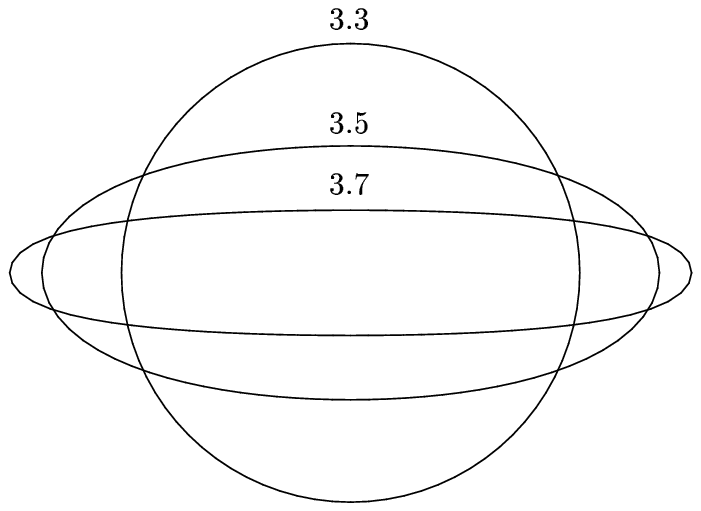}}} 
         \hspace{.5in}
  \vbox{\hbox{\includegraphics[height=1.5in]{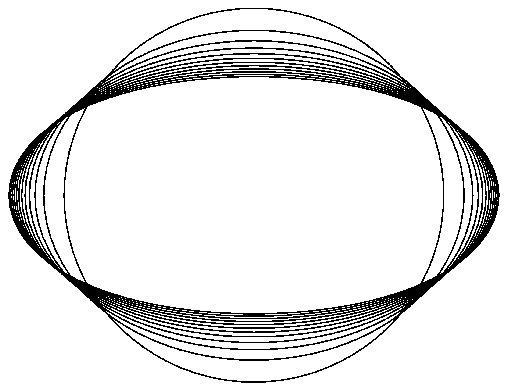}}}
  \hfill} 
\caption[A collection of maximizers]{A collection of curves of length
$2\pi$ which maximize average chord length to the $p$-th power for
various values of $p$. The curves on the left are labelled with the
corresponding values of $p$. The curves on the right represent values
of $p$ from $3.462$ to $3.484$ in increments of $0.002$. These curves
are numerical approximations of the true maximizers computed with
Brakke's {\em Evolver}. }\label{fig:curvefig}
\end{figure}

In any event, it is interesting to consider how the shape of the
maximizers changes as we vary~$p$. Since the limit of~$L^p$ norms
as~$p\to\infty$ is the supremum norm, we have
\begin{equation*}
\lim_{p\to \infty}
A_p[c]
=
\sup_{s,t} |c(t)-c(s)| \le \pi
\end{equation*}
with equality if and only if~$c$ double covers a line segment of
length~$\pi$.  So the $c_p$ form a family of convex curves converging
to the double-covered segment as $p \to \infty$, and to the circle as
$p \to 2$.  To illuminate this process, we numerically computed
maximizers of $A_p$ for values of $p$ between~$2$ and~$4$ using
Brakke's Evolver~\cite{brakke:evolver}. Figure~\ref{fig:curvefig}
shows some of the~$c_p$.

Since the double-covered segment has greater average $p$-th power
chord length than the circle for $p > 3.5721$, there
must be some critical value $p^*$ of $p$ between $2$ and $3.5721$
where ``the symmetry breaks'', and circles are no longer maximizers
for $A_p$.

To find an approximate value for $p^*$, we computed the ratio $r(p)$
of the widest and narrowest projections of each of our computed
maximizers for $p$ between $2$ and $4$.  Since all these curves are
convex, a value close to unity indicates a curve close to a circle.

\figh{widthfig}{2}{This figure shows two plots of the ratio~$r(p)$ of
the widest and narrowest projections of the computed maximizers of
average chord length to the $p$-th power for values of~$p$ between~$1$ and~$4$.}

As Figure~\ref{fig:widthfig} shows, by this measure the computed
minimizers are numerically very close to circles for $2 \leq p \leq
3.45$. To check this conclusion, we fit each minimizer to an ellipse
using a least-squares procedure.  Figure~\ref{fig:fitfig} shows the
results of these computations.

\figh{fitfig}{2}{The base-$10$ logarithm of the error~$e(p)$ in a
least-squares fit of the computed maximizer for average chord length to 
the $p$-th power to an ellipse, plotted against $p$.}

To give a sense of the accuracy of our computations, this graph
includes some computed minimizers for $p$ between $1$ and $2$, for
which we have proved that the unique minimizer is the 
circle. We also computed the eccentricities of each of the best-fit
ellipses. 

A conservative reading of all this data supports the surprising
conjecture that $p^*$ is at least $3.3$. Further, we note that for $p
> p^*$, the maximizing curves do not seem to be ellipses, as one might
have conjectured by looking at Theorem~\ref{kostya}.

{\bf Acknowledgments.} We thank Kostya Oskolkov for pointing out that the
complex form of Fourier series would simplify the first proof of
Theorem~\ref{kostya}.  We also owe a debt to the bibliographic
notes in the wonderful book of Santal{\'o}~\cite{santalo:book} for the
reference to \Luko's paper~\cite{Gabor}.

%\bibliographystyle{abbrv}
%\bibliography{references}

%\bibliographystyle{amsplain}
%\bibliography{/math/faculty/howard/tex/inputs/HowRefs,references}

\begin{thebibliography}{10}

\bibitem{brakke:evolver}
K. Brakke, \emph{The Surface Evolver}, Experimental Math. \textbf{1} (1992),
no.~2, 141--165.

\bibitem{Burago-Zal:geometric-ineq}
J.~D. Burago and V.~A. Zalgaller, \emph{Geometric inequalities}, Grundlehren,
  vol. 285, Springer, Berlin, 1980.

\bibitem{chern:book}
S.~S. Chern, \emph{Curves and surfaces in {E}uclidean space}, Studies in Global
  Geometry and Analysis, Math. Assoc. Amer. (distributed by Prentice-Hall,
  Englewood Cliffs, N.J.), 1967, pp.~16--56. \MR{35 \#3610}

\bibitem{freedman&he&wang:mobius}
M.~H. Freedman, Z.-X. He, and Z.~Wang, \emph{M\"obius energy of knots and
  unknots}, Ann. of Math. (2) \textbf{139} (1994), no.~1, 1--50. \MR{94j:58038}

\bibitem{Ghomi-Howard;unfold}
M.~Ghomi and R.~Howard, \emph{Convex unfoldings of space curves}, Preprint.

\bibitem{gromov:filling}
M.~Gromov, \emph{Filling {R}iemannian manifolds}, J. Differential Geom.
  \textbf{18} (1983), no.~1, 1--147. \MR{85h:53029}

\bibitem{Gromov:metric}
\bysame, \emph{Metric structures for {R}iemannian and non-{R}iemannian spaces},
  Birkh\"auser Boston Inc., Boston, MA, 1999, Based on the 1981 French
  original, With appendices by M.\ Katz, P.\ Pansu and S.\ Semmes, Translated
  from the French by Sean Michael Bates. \MR{2000d:53065}

\bibitem{hurwitz:iso}
A.~Hurwitz, \emph{Sur le probl\'eme des isop\`erim\'etres}, C. R. Acad. Sci.
  Paris \textbf{132} (1901), 401--403, Reprinted in
  \cite[pp.~490--491]{Hurwitz:work1}.

\bibitem{Hurwitz:work1}
\bysame, \emph{Mathematische {W}erke. {B}d. {I}: {F}unktionentheorie},
  Birkh\"auser Verlag, Basel, 1962, Herausgegeben von der Abteilung f\"ur
  Mathematik und Physik der Eidgen\"ossischen Technischen Hochschule in
  Z\"urich. \MR{27 \#4723a}

\bibitem{kusner&sullivan:distortion}
R.~B. Kusner and J.~M. Sullivan, \emph{On distortion and thickness of knots},
  Topology and geometry in polymer science (Minneapolis, MN, 1996), Springer,
  New York, 1998, pp.~67--78. \MR{99i:57019}

\bibitem{Gabor}
G.~\Luko, \emph{On the mean length of the chords of a closed curve}, Israel
  J. Math. \textbf{4} (1966), 23--32. \MR{34 \#681}

\bibitem{ohara:energy1}
J.~O'Hara, \emph{Energy of a knot}, Topology \textbf{30} (1991), no.~2,
  241--247. \MR{92c:58017}

\bibitem{ohara:energy2}
\bysame, \emph{Energy Functionals of Knots}, In Proc. Topology Conference Hawaii,
  Karl Dovermann, ed. World Scientific, New York, 1992, pp.~201--214.

\bibitem{ohara:energy2a}
\bysame, \emph{Family of energy functionals of knots}, Topology Appl.
  \textbf{48} (1992), no.~2, 147--161. \MR{94h:58064}

\bibitem{ohara:energy3}
\bysame, \emph{Energy functionals of knots. {I}{I}}, Topology Appl. \textbf{56}
  (1994), no.~1, 45--61. \MR{94m:58028}

\bibitem{Royden:analysis}
H.~L. Royden, \emph{Real analysis}, third ed., Macmillan Publishing Company,
  New York, 1988. \MR{90g:00004}

\bibitem{sallee}
G.~T. Sallee, \emph{Stretching chords of space curves}, Geometriae Dedicata \textbf{2} (1973),
311--315. \MR{49 \#1334 53A05}

\bibitem{santalo:book}
L.~A. Santal{\'o}, \emph{Integral geometry and geometric probability},
  Addison-Wesley Publishing Co., Reading, Mass.-London-Amsterdam, 1976, With a
  foreword by Mark Kac, Encyclopedia of Mathematics and its Applications, Vol.
  1. \MR{55 \#6340}

\bibitem{schneider:book}
R. Schneider, \emph{Convex bodies: the {B}runn-{M}inkowski theory}, Cambridge
  University Press, Cambridge, 1993. \MR{94d:52007}

\bibitem{Spivak:comp2}
M.~Spivak, \emph{A comprehensive introduction to differential geometry}, 2 ed.,
  vol.~2, Publish or Perish Inc., Berkeley, 1979.

\bibitem{Spivak:comp5}
\bysame, \emph{A comprehensive introduction to differential geometry}, 2 ed.,
  vol.~5, Publish or Perish Inc., Berkeley, 1979.

\end{thebibliography}

\providecommand{\bysame}{\leavevmode\hbox to3em{\hrulefill}\thinspace}
\providecommand{\MR}{\relax\ifhmode\unskip\space\fi MR }
\providecommand{\MRhref}[2]{%
  \href{http://www.ams.org/mathscinet-getitem?mr=#1}{#2}
}
\providecommand{\href}[2]{#2}

\end{document}